\documentclass{amsart}
\usepackage{amssymb, latexsym}

\newcommand{\XREF}{\label}

\newcommand{\C}{{\mathbb C}}

\newcommand{\CC}{{\overline \C}}

\newcommand{\A}{{\mathcal A}}
\newcommand{\downto}{{\searrow}}
\newcommand{\F}{{\mathcal F}}

\newcommand{\nbhd}{{neighborhood }}
\newcommand{\tri}{{\bigtriangleup}}

  \newtheorem{corollary}{Corollary}
\newtheorem{lemma}{Lemma}

\newtheorem{theorem}{Theorem}

\theoremstyle{definition}
\newtheorem{definition}{Definition}
\theoremstyle{remark}
\newtheorem{remark}{Remark}

\newtheorem{example}{Example}

\numberwithin{equation}{section}
\numberwithin{theorem}{section}
\numberwithin{definition}{section}
\numberwithin{remark}{section}
\numberwithin{example}{section}
\numberwithin{lemma}{section}
\numberwithin{property}{section}
\numberwithin{proposition}{section}
\numberwithin{claim}{section}
\numberwithin{othertheorem}{section}
\numberwithin{conj}{section}
\numberwithin{corollary}{section}

\begin{document}

\title{Uniformly perfect analytic and conformal attractor sets}
\author{Rich Stankewitz}
\address{Department of Mathematics,
         Texas A\&M University,
         College Station, TX 77843,
    richs@math.tamu.edu}

\curraddr{School of Science\\
    Penn State University at Erie\\
    Station Road\\
    Erie, PA 16563\\
    rls42@psu.edu}

\begin{abstract}
Conditions are given which imply that analytic iterated function
systems (IFS's) in the complex plane $\C$ have uniformly perfect
attractor sets.  In particular, it is shown that the attractor set of a
finitely generated conformal IFS is uniformly perfect when it contains
two or more points.  Also, an example of a finitely generated
analytic attractor set which is not uniformly perfect is given.
\end{abstract}

\maketitle

\section{Introduction and results}
Consider an iterated function system
(IFS) $G=\langle g_i:i \in I
\rangle$, the set of all finite compositions of non-constant
generating maps $\{g_i:i \in I\}$ for some index set $I$, where each
function maps the open connected set $U
\subset \C$ into a compact set $K \subset U$ such that there exists
$0<s<1$ and a metric $d$ on $K$ where $d(g_i(z),g_i(w))\leq s d(z,w)$
for all $z, w \in K$ and all $i \in I$.  Thus this system is uniformly
contracting on the metric space $(K,d)$.  We say that $G$ is an IFS on
$(U,K)$.  We define the attractor set $A$ by $A=A(G)=\overline{A'}$,
the closure of $A'$ in the Euclidean topology, where
$A'=A'(G)=\{z:\text{there exists } g
\in G \text{ such that } g(z)=z
\}$ is the set of (attracting) fixed points of $G$.  We
will suppress the dependence on $G$ when there is no chance for
confusion.  We define an IFS and its corresponding attractor set
to be analytic (respectively conformal) if all the maps are analytic
(respectively conformal) on $U$.  In an analytic IFS, we note that we
may, and will, take the metric $d$ to be the hyperbolic metric on $U$
(see section 2).

We note that when the IFS $G=\langle g_1, \dots, g_N \rangle$, we have
that $A$ is the unique compact set such that
\begin{equation}\label{invariance}
A = \cup_{i=1}^N g_i(A)
\end{equation}
(see~\cite{Hu}, p.~724).  We also point out that in~\cite{MU1,MU2} the
limit set $J$ of a conformal IFS is defined in such a way so that $J
\supset A'$ and $\overline{J}=A$.  In~\cite{MU1,MU2} conformal
IFS's are defined differently, and non equivalently, and so the reader
should observe that the definition above does not match exactly the
definitions that one may find in the literature.

In~\cite{RS3} certain conformal attractor sets have been shown to be
uniformly perfect, that is, these sets are uniformly thick near each
of their points (see section 2), when the generating maps are M\"obius.
We note that in~\cite{RS3}, the result in Corollary 3.3 is only stated
for attractor sets generated by linear maps, but clearly more general
cases are handled by Theorem 3.1.  In this reference critical use was
made of the fact that the generating maps were globally defined (and
had globally defined inverses).  Such restrictions are not needed in
the following theorem.

\begin{theorem}\label{better}
Let $G=\langle g_i:i \in I \rangle$ be an analytic IFS on $(U,K)$
where there exist $0<\delta<diam(A)$ and $C>0$ such that we have the
following:

(I) if $a \in A$ and $i \in I$, then $g_i$ is one-to-one on
$\tri(a,\delta)$,

(II) if $a \in A$ and $g_i(a)=a'$, then the branch $h_i$ of
$g_i^{-1}$ such that $h_i(a')=a$ is defined on $\tri(a',\delta),$

(III) if $a \in A$ and $g_i(a)=a'$, then the branch $h_i$ of
$g_i^{-1}$ such that $h_i(a')=a$ satisfies

$|h_i(z)-h_i(a')| \leq C|z-a'|$ for all $z \in
\overline{\tri(a',\delta/10)}$.

Then, if the attractor set $A$ has infinitely many points,
then $A$ is uniformly perfect.
\end{theorem}
Above and throughout we use the following notation.
Let $q$ be a metric.
For a set $F\subset
\C$, let $diam_q(F)=\sup\{q(z,w):z, w \in F\}$ and
$F^q_\epsilon = \{z:dist_q(z,F) <\epsilon\}$ where
$dist_q(z,F)=\inf\{q(z,w):w \in F\}$.  Also let
$\tri_q(w,r)=\{z:q(z,w)<r\}$ and $C_q(w,r)=\{z:q(z,w)=r\}.$ If no
metric is noted, then it is assumed that the metric is the
Euclidean metric.

\begin{remark}\label{gis1-1}
Suppose $\delta$ satisfies (I) in Theorem~\ref{better}.
Then there exists (see proof at the end of section 2)
$\delta'>0$ such that each $g \in G$ is one-to-one on
$\tri(a,\delta')$.  Hence in the proof of Theorem~\ref{better} we will
replace $\delta$ by $\delta'$ and assume that each $g \in G$ is
one-to-one on $\tri(a,\delta)$ when $a \in A$.
\end{remark}

\begin{remark}\label{replace}
Suppose $\delta$ satisfies (I) in Theorem~\ref{better}.  If there
exists $\epsilon>0$ such that for all $a \in A$ and $i \in I$ we have
 $g_i(\tri(a,\delta))
\supset \tri(g_i(a),\epsilon)$, then we see that the branch $h_i$ of
$g_i^{-1}$ such that $h_i(g_i(a))=a$ is defined on
$\tri(g_i(a),\epsilon)$.  Hence we may replace $\delta$ by
$\min\{\delta,\epsilon\}$ to satisfy both (I) and (II).
\end{remark}

\begin{corollary}\label{cor3}
Let $G=\langle g_i:i \in I \rangle$ be an analytic IFS on $(U,K)$ such
that there exists $\eta>0$ where $|g_i'(a)| \geq \eta$ for all $a \in
A$ and all $i \in I$.  If $A$  has infinitely many points,
then $A$ is uniformly perfect.
\end{corollary}

\begin{remark}
We show in Lemma~\ref{equivhyp} that the hypotheses of
Theorem~\ref{better} and Corollary~\ref{cor3} are equivalent.
\end{remark}


\begin{corollary}\label{main}
Let $G=\langle g_i:i \in I \rangle$ be a conformal IFS on $(U,K)$ such
that there exist $\eta>0$ where $|g_i'(a)| \geq \eta$ for all $a \in
A$ and all $i \in I$.
If $A$ contains more than one point, then $A$ is
uniformly perfect.
\end{corollary}

\begin{corollary}\label{cor}
Let $G=\langle g_1, \dots, g_N \rangle$ be a conformal IFS on $(U,K)$.
If $A$ contains more than one point, then $A$ is uniformly perfect.
\end{corollary}

Uniformly perfect sets, which are defined in section~2, were
introduced by A. F. Beardon and Ch. Pommerenke in 1978 in~\cite{BP}.
Such sets cannot be separated by annuli that are too ``fat'' (large
ratio of outer to inner radius) and thus near each of its points the
set is uniformly ``thick'' where this thickness is independent of
scaling.  See~\cite{BP,P,P2} for many interesting equivalent
definitions.  Such sets are known to be regular for the Dirichlet
problem (see~\cite{P}, p.~193) and also have positive Hausdorff
dimension (see~\cite{JV}, p.~523).  See~\cite{MU1,MU2} for a
discussion on the Hausdorff dimension, packing dimension, and other
properties of limit sets (attractors) for conformal IFS's.

This paper is organized as follows. Section~2 contains basic lemmas
and definitions.  Section~3 contains the proofs of
Theorem~\ref{better} and
Corollaries~\ref{cor3},~\ref{main} and ~\ref{cor}.
Section~4 contains examples of the following: a uniformly perfect
conformal attractor set that is generated by maps of the form $z^2+c$;
a uniformly perfect analytic attractor set that is generated by maps
of the form $z^2+c$; a non uniformly perfect attractor set which is
generated by an infinite number of conformal maps; a (non-conformal)
IFS generated by two analytic maps whose attractor set is perfect
though not uniformly perfect, yet the generating maps are one-to-one
on $A$; analytic attractor sets whose cardinality is any given integer
$n$.

The author would like to Harold Boas, Marshall
Whittlesey, Mohammed Ziane, and especially Aimo Hinkkanen for
productive conversations regarding this subject.

\section{Definitions and basic facts}
Given an open set $U \subset \C$ such that there exists a non constant
analytic map $g$ with $g(U)$ contained in a compact set $K \subset U$,
then $U$ is hyperbolic, that is, $\C \setminus U$ contains more than
one point.  This follows since the image of the plane or the punctured
plane under an analytic map is always dense in the plane and hence
could not be contained in a compact subset of $U$.  Hence there exists
a hyperbolic metric on $U$ (see~\cite{CG}, p.~12).

\begin{lemma}\label{hypmetric}
If the analytic function $g$ maps an open connected set $U \subset
\C$ into a compact set $K \subset U$, then there exists $0<s<1$,
which depends on $U$ and $K$ only, such that $d(g(z),g(w))\leq s
d(z,w)$ for all $z, w \in K$ where $d$ is the
hyperbolic metric defined on $U$.
\end{lemma}

\begin{proof}
If no such $s$ exists then there exist sequences $z_n, w_n \in K$ and
analytic maps $g_n$ such that $g_n(U) \subset K$ where
$d(g_n(z_n),g_n(w_n))/d(z_n,w_n) \to 1$.  By compactness and normality
we may assume that $z_n \to z, w_n \to w$ and $g_n \to g$ uniformly on
$K$ where $g(U) \subset K$.  If $z \neq w$ then one can show
$d(g_n(z_n),g_n(w_n))/d(z_n,w_n)
\to d(g(z),g(w))/d(z,w) <1$ (see~\cite{CG}, p.~12).  If $z=w$, then then
$d(g_n(z_n),g_n(w_n))/d(z_n,w_n)
 \to |g'(z)|\lambda(g(z))/\lambda(z)<1$ where
$\lambda(z)$ denotes the density of the hyperbolic metric of $U$ at the
point $z$ of $U$.  Thus we have a contradiction and so the lemma follows.
 \end{proof}



\begin{lemma}\label{forinv}
Let $G=\langle g_i:i \in I \rangle$ be an IFS on $(U,K)$.  Then for
any $g \in G$, we have $g(A) \subset A$.
\end{lemma}

\begin{proof}
Suppose $a \in A'(G)$ and so $f(a)=a$ for some $f \in G$ with, say,
$f=g_{i_1} \circ g_{i_2} \circ \dots \circ g_{i_k}$ where each
$i_j \in I$.  Say $g=g_{i_1'} \circ g_{i_2'} \circ
\dots \circ g_{i_n'}$ where each $i_j' \in I$.
Letting $\tilde{A}= A(\langle g_{i_1}, g_{i_2},\dots, g_{i_k}, g_{i_1'},
g_{i_2'}, \dots, g_{i_n'} \rangle)$ we see by~\eqref{invariance} that $g(a)
\in g(\tilde{A}) \subset \tilde{A} \subset A(G).$ Hence $g(A'(G))
\subset A(G)$ and the lemma follows from the continuity of $g$.
\end{proof}



\begin{lemma}\label{gi1-1}
Let $\eta>0, r>0, M>0, a \in \C$ and let $\F$ denote the family of maps
$\{f:\tri(a,r) \to
\tri(0,M) \text{ such that } f \text{ is analytic and }|f'(a)| \geq \eta\}$.
Then there exists $\rho>0$ such that each $f
\in \F$ is one-to-one on $\tri(a,\rho)$.
\end{lemma}

\begin{remark}
Any $\rho$ such that $M
\rho(2r-\rho)/(r(r-\rho)^2) <\eta$ will
suffice will satisfy the conclusion of Lemma~\ref{gi1-1}.
\end{remark}

\begin{proof}
Fix $\rho$ such that $M \rho(2r-\rho)/(r(r-\rho)^2) <\eta$.  For
$|z-a| < r$ we may obtain by estimates on the Cauchy integral formula
that $|f''(z)| \le 2 M r/(r-|z-a|)^3$.  Hence for $|z-a|
\le \rho$ we have
$|f'(z)-f'(a)|=|\int_a^z f''(s) \,\,ds| \le 2 M r \int_0^{|z-a|}
(r-t)^{-3} \,\,dt = Mr[(r-|z-a|)^{-2} - r^{-2}] \le M
\rho(2r-\rho)/(r(r-\rho)^2) <\eta \le |f'(a)|$ where the integral is
over the straight line path.  This suffices to show that $f$ is
one-to-one on $\tri(a,\rho)$ (see~\cite{Co}, p.~293).
\end{proof}

\begin{lemma}\label{equivhyp}
Let $G=\langle g_i:i \in I \rangle$ be an analytic IFS on $(U,K).$
Then there exist $0<\delta<diam(A)$ and $C>0$ such that conditions (I)-(III)
of Theorem~\ref{better} hold if and only if there exists $\eta>0$ such
that we have $|g_i'(a)| \geq \eta$ for all $a \in A$ and all $i \in
I$.
\end{lemma}

\begin{proof}
Let $0<\delta<diam(A)$ and $C>0$ be such that conditions (I)-(III)
hold.  If there does not exist $\eta>0$ such that we have $|g_i'(a)|
\geq \eta$ for all $a \in A$ and all $i \in I$, then there exist $a_n
\in A$ and $g_n \in \{g_i:i \in I\}$ such that $g_n'(a_n)\to 0$.  Let
$a_n'=g_n(a_n)$ and let $h_n$ be the branch of
the inverse of $g_n$ such that $h_n(a_n')=a_n$.  Conditions (II) and (III)
imply, in particular, that $h_n(\tri(a_n',\delta/10)) \subset
\tri(a_n,C\delta/10)$.  But by the Koebe 1/4 theorem we have
$h_n(\tri(a_n',\delta/10)) \supset
\tri(a_n,|h_n'(a_n')|\delta/40)$ which gives a contradiction for
large $n$ since $|h_n'(a_n')|=1/|g_n'(a_n)|\to \infty.$

Suppose $\eta>0$ is such that $|g_i'(a)| \geq \eta$ for all $a
\in A$ and all $i \in I$.  Note that since $K$ is compact
 there exists $M>0$ such that each $g_i$ maps $U$ into some
 $\tri(0,M)$.  Let $r>0$ be such that $A_r \subset U$.  By
 Lemma~\ref{gi1-1} setting $\delta=\rho$ we see that condition (I) is
 satisfied.  By the Koebe 1/4 theorem each $g_i(\tri(a,\delta))
 \supset
\tri(g_i(a),\delta\eta/4)$ and so (see Remark~\ref{replace}) we may
replace $\delta$ by $\min\{\delta,\delta\eta/4\}$
to satisfy both (I) and (II).
Since the branch $h_i$ of $g_i^{-1}$ such that $h_i(g_i(a))=a$ is
defined on $\tri(g_i(a),\delta)$ we see that
$|h_i'(g_i(a))|=1/|g_i'(a)| \leq 1/\eta$ and so by Lemma~\ref{Lip} we
see that condition (III) is satisfied with $C=100/(81 \eta)$.
\end{proof}

\begin{lemma}\label{Aperfect}
Let $G=\langle g_i:i \in I \rangle$ be an IFS such that for each $g
\in G$ there exists a positive integer $N_g$ such that for each $a \in A$
the set $g^{-1}(\{a\})$ has at most $N_g$ elements.  If $A$ has
infinitely many elements, then $A$ is perfect and hence, by the Baire
category theorem is uncountable.
\end{lemma}

\begin{proof}
If $a\in A=\overline{A'}$ were isolated it would have to be an
attracting fixed point for some element $f \in G$.  Suppose $\tri(a,r)
\cap A = \{a\}$.  For some positive integer $n$, we have $f^n(A) \subset
\tri(a,r)$
since $a \in f^n(A)$ and $diam_d(f^n(A)) \to 0$ as $n
\to \infty$.  By Lemma~\ref{forinv} we see that $f^n(A)=\{a\}$ and
this contradicts the hypothesis on the integer $N_{f^n}$.
\end{proof}

\begin{remark}\label{idthm}
If $G$ is an analytic IFS, then the proof of Lemma~\ref{Aperfect} and
the identity theorem for non-constant analytic maps imply
that the attractor set $A$ is perfect when $A$ has infinitely many
elements.
\end{remark}

\begin{lemma}\label{backmaps}
Let $G=\langle g_i:i \in I \rangle$ be an analytic IFS such that each $g_i'$
never takes the value zero on $A$.  If $a \in
A'$, then there exists two sequences of maps $F_1=g_{i_1},
F_2=g_{i_1} \circ g_{i_2}, \dots, F_m=g_{i_1} \circ g_{i_2} \circ
\dots \circ g_{i_m}, \dots$ where each $i_j \in I$ and
$H_1=h_{i_1}, H_2=h_{i_2} \circ h_{i_1},\dots,
H_m=h_{i_{m}} \circ \dots \circ h_{i_1}, \dots$ where each $h_{i_j}$ is
a branch of $g_{i_j}^{-1}$, such that for each positive integer $m$ we have
$(H_m \circ F_m)(z)=z$ (on a small neighborhood of $H_m(a)$)
and $H_m(a) \in A.$
\end{lemma}

\begin{proof}
Suppose $g(a)=a$ for $g \in G$ and so $g=g_{i_1} \circ g_{i_2} \circ
\dots \circ g_{i_k}$ for some choice of $i_j \in I$.  Let $h_{i_1}$ be the
branch of $g_{i_1}^{-1}$ such that $h_{i_1}(a)= (g_{i_2} \circ \dots
\circ g_{i_k})(a) \in A$ by Lemma~\ref{forinv}.
Similarly choose $h_{i_2},\dots,h_{i_k}$ to be branches of
$g_{i_2}^{-1},\dots,g_{i_k}^{-1}$, respectively, so that
$h_{i_2}(h_{i_1}(a))= (g_{i_3} \circ \dots \circ g_{i_k})(a) \in A,
\dots,   h_{i_j}((h_{i_{j-1}}\circ \dots \circ h_{i_1})(a))=(g_{i_{j+1}} \circ \dots \circ g_{i_k})(a) \in A, \dots,  h_{i_k}((h_{i_{k-1}}\circ \dots \circ h_{i_1})(a))=a \in A.$

We define $F_1=g_{i_1}, F_2=g_{i_1}
\circ g_{i_2}, \dots, F_k=g_{i_1} \circ g_{i_2} \circ \dots \circ
g_{i_k}, F_{k+1}=g_{i_1} \circ g_{i_2} \circ \dots \circ
g_{i_k}\circ g_{i_1}, \dots$.  For $F_m=g_{i_1} \circ g_{i_2} \circ
\dots \circ g_{i_m}$ we define $H_m=h_{i_{m}} \circ \dots \circ h_{i_1}$
and we are done.
\end{proof}

\begin{definition} \XREF{sepann}
A true annulus $\A=Ann(w;r,R)=\{z:r < |z-w| < R\}$ is said to
separate a set $F$ if $F$ intersects both components of $\CC \setminus
\A$ and $F \cap \A = \emptyset$.
\end{definition}

\begin{definition} \XREF{updef}[\cite{P}, p.~192]
A compact subset $F \subset \C$ is uniformly perfect if
there exists $M>1$ such that whenever $w \in F,\, 0<R<diam(F),\, R/r > M,$
then $Ann(w;r,R) \cap F \neq \emptyset$.
\end{definition}

\begin{lemma}\label{Rntozero}[\cite{P2}, p. 315]
Let $F$ be a perfect compact subset of $\C$ and let
$A_n=Ann(z_n;r_n,R_n)$ be annuli which separate $F$, have $z_n \in
F$ and $R_n/r_n
\to \infty$.  Then $R_n \to 0$.
\end{lemma}

\begin{proof}
Since $F$ is compact we may assume that $z_n \to w \in F$.  We also
see that $R_n \leq diam(F)$ and hence $r_n \to 0$.  If there exists
$\delta >0$ such that $R_n > \delta$ along some subsequence of $R_n$
then one can show that $\tri(w,\delta) \cap F=\{w\}$ since
$\tri(w,\delta) \setminus \{w\}$ is contained in the ``limit'' of this
subsequence of $A_n$.  This contradicts the hypothesis that $F$ is
perfect.
\end{proof}

\begin{lemma}\label{domain}
Let $\delta>0$ and let $g$ be analytic and univalent on $\tri(0,\delta)$ with
$g(0)=0$ such that the there exists $z$ with $|z| \leq \delta/10$ and
$|g(z)|=R$. Then $g(\tri(0,\delta)) \supset \tri(0, 2 R)$ and hence
$g^{-1}$ is defined on $\tri(0, 2 R)$.
\end{lemma}

\begin{proof}
Defining the univalent map $K(z)=\frac{1}{\delta g'(0)} g(\delta z)$
on $\tri(0,1)$ we have $K(0)=0$ and $K'(0)=1$.  By hypothesis there
exists a point $z'=z/\delta$ such that $|z'| \leq 1/10$ and $|K(z')|=
\frac{R}{\delta |g'(0)|}$.  By the Koebe distortion theorem and the maximum
modulus theorem (see~\cite{CG}, p.~3)
we see that for $|z|\leq 1/10$ we have
$$|K(z)| \leq \frac{1/10}{(1-1/10)^2}$$
and so $$\frac{R}{\delta |g'(0)|} \leq \frac{10}{81}.$$

By the Koebe 1/4 theorem (see~\cite{CG}, p.~2) we have
$K(\tri(0,1)) \supset \tri(0,1/4)$ and so we
conclude that $g(\tri(0,\delta)) \supset
\tri(0,\delta |g'(0)|/4) \supset \tri(0,2R)$.
\end{proof}

\begin{lemma}\label{Lip}
Let $\eta>0$ and let $g$ be analytic and univalent on $\tri(a,\delta)$ with
 $|g'(a)| \leq 1/\eta$. If $|z-a| \leq \delta/10$, then
$|g(z)-g(a)| \leq 100 |z-a|/(81 \eta).$
\end{lemma}

\begin{proof}
This follows from the Koebe distortion theorem in much the same manner
as the proof of Lemma~\ref{domain} and therefore we omit the details.
\end{proof}

\begin{lemma}\label{annlemma}
Let $g$ be univalent on $\tri(0,2R)$ with $g(0)=0$.  If $R>9r$, then
$g(Ann(0;r,R)) \supset Ann(0;r',R')$ for some $R'>r'>0$ such that
$R'/r' \geq R/(9r)$.
\end{lemma}

\begin{proof}
The univalent map $h(z)= \frac{g(2Rz)}{2R g'(0)}$ defined in $\tri(0,1)$
is such that
$h(0)=0$ and $h'(0)=1$.  By the Koebe distortion theorem we have $$|h(z)|
\leq \frac{r/(2R)}{[1-(r/(2R))]^2}  \text{ on } |z|=r/(2R) \qquad \text{ and
}\qquad |h(z)| \geq \frac{1/2}{[1+(1/2)]^2} \text{ on }|z|=1/2.$$
So $$|g(z)| \leq \frac{r|g'(0)|}{[1-(r/(2R))]^2} =r'\text{ on }
|z|=r \qquad \text{ and }\qquad |g(z)| \geq \frac{R|g'(0)|}{[1+(1/2)]^2}=R'
\text{ on }|z|=R.$$

Hence $g(Ann(0;r,R)) \supset Ann(0;r',R')$ where $\frac{R'}{r'} =
\frac{R [1-(r/(2R))]^2}{r [1+(1/2)]^2}\geq R/(9r)$ since $r<R$.
\end{proof}

\begin{lemma}\label{metrics}
Let $\rho>0$.  If $p$ and $q$ are metrics that induce the same
topology on a compact topological space $F$, then there exists
$r>0$ such that $\tri_p(a,r) \subset \tri_q(a,\rho)$ for all $a \in F$.
\end{lemma}

\begin{proof}
Consider the open cover $\{\tri_q(a,\rho/2):a \in F\}$ of $F$.  By the
Lebesgue number lemma (see~\cite{Mu}, p.~179) applied to the compact
metric space $(F,p)$, there exists $R>0$ such that for each $a \in F$
there exists $a' \in F$ such that $\tri_p(a,R/3) \subset
\tri_q(a',\rho/2) \subset \tri_q(a,\rho)$.  Hence we see the result
follows with $r=R/3$.
\end{proof}

\begin{proof}[Proof of Remark~\ref{gis1-1}]
Let $d$ be the hyperbolic metric on $U$ and recall that the hyperbolic
and Euclidean metrics induce the same topology on $U$.  Since $A$ is
compact one may apply Lemma~\ref{metrics} twice to obtain $\delta',
r>0$ such that $\tri(a,\delta') \subset \tri_d(a,r) \subset
\tri(a,\delta)$ for all $a \in A$.
Hence for any $f \in G$ we have that
$f(\tri(a,\delta')) \subset f(\tri_d(a,r)) \subset
\tri_d(f(a),r) \subset
\tri(f(a),\delta)$.  Note that
we used the fact that $f(A) \subset A$.

Let $g \in G$ and say it can be written as $g=g_{i_k} \circ \dots
\circ g_{i_2} \circ g_{i_1}$.  Note that $g_{i_1}$ is one-to-one on
$\tri(a,\delta')$.  Since from above $g_{i_1}(\tri(a,\delta')) \subset
\tri(g_{i_1}(a),\delta)$ we use (property (I)) the fact that $g_{i_2}$
is one-to-one on $\tri(g_{i_1}(a),\delta)$ to conclude that
$g_{i_2} \circ g_{i_1}$ is one-to-one on $\tri(a,\delta')$ and again
from above $(g_{i_2} \circ g_{i_1})(\tri(a,\delta')) \subset
\tri((g_{i_2} \circ g_{i_1})(a),\delta)$ .  We then
proceed inductively to show that each map $g_{i_m} \circ \dots
\circ g_{i_2} \circ g_{i_1}$ including $g$ is one-to-one on $\tri(a,\delta')$.
\end{proof}

\section{Proof of the Theorem and Corollaries}

\begin{proof}[Proof of Theorem~\ref{better}]



Let the hypotheses of the theorem be satisfied.  Recall that $d$ is
the hyperbolic metric on $U$.  Since $K$ is a compact subset of $U$
there exists $\epsilon >0$ such that $\overline{K^d_\epsilon} \subset
U$.  Let $\tilde{s}<1$ be the uniform contraction coefficient on
$\overline{K^d_\epsilon}$ for the generating maps as given by
Lemma~\ref{hypmetric}.  Since we may shrink $\delta$ and still satisfy
conditions (I)-(III) of the theorem we will assume (see
Lemma~\ref{metrics}) that for each $a \in A$ we have
$\tri(a,\delta/(10C)) \subset \tri_d(a,\epsilon)$.

Suppose $A$ contains infinitely many points, but is not uniformly
perfect.  Hence there exist
annuli $A_n=Ann(a_n;r_n,R_n)$ that separate $A$ where $R_n/r_n
\to \infty$ and each $a_n \in A$.  Since $A'$ is dense in $A$ a simple
geometric argument allows one to assume that each $a_n \in A'$.  Since
$A$ is perfect (see Remark~\ref{idthm}) we must have both $r_n$ and
$R_n$ tending to zero by Lemma~\ref{Rntozero}. We assume all $R_n <
\delta/(10C)$ and $R_n>9r_n$.

The method of proof will be to show that the $A_n$ can be expanded by
certain locally defined inverses of maps in $G$ given in
Lemma~\ref{backmaps}, such that the expanded (conformal) annuli
contain true annuli of large ratio of outer to inner radius.  The
Koebe 1/4 theorem (Lemma~\ref{domain}) and the hypothesis on $\delta$
is used to show that the chosen inverse maps are well defined on
$\tri(a_n,2R_n)$.  The deformation that occurs since the maps are not
linear (similitudes) is controlled by the use of the Koebe distortion
theorem (see Lemma~\ref{annlemma}).  The expanded true annuli are
constructed so that the outer radii do not tend to zero, yet these
annuli do separate $A$.  This contradicts the fact that $A$ is perfect
(see Lemma~\ref{Rntozero}).  We now carry out this program.

Since $a_n \in A'$ we may construct sequences of maps $F_m$ and $H_m$ for
$m=1, 2, \dots$ as given by Lemma~\ref{backmaps}.  Note that the maps
$F_m$ and $H_m$ depend on $n$.
There exists $m$ such that
\begin{equation}\label{min}
F_m(\tri(H_m(a_n),\delta/(10C))) \subset F_m(\tri_d(H_m(a_n),
\epsilon)) \subset \tri_d(a_n, \tilde{s}^m
\epsilon)
\subset\tri(a_n,R_n).
\end{equation}
The first inclusion follows from the choice of $\delta$ as stated at
the beginning of the proof and the last inclusion follows from the
fact that the Euclidean and the hyperbolic metrics induce the
same topology on $U$.

We define $f_n=F_{m^*}$ where $m^*$ is the smallest integer $m$ such
that $F_m(\tri(H_m(a_n),\delta/(10C))) \subset \tri(a_n,R_n)$
and define its ``inverse'' by $k_n=H_{m^*}$.  Therefore
$F_{{m^*}-1}(\tri(H_{{m^*}-1}(a_n),\delta/(10C))) \cap C(a_n,R_n) \neq
\emptyset.$
Choose $z \in \tri(H_{{m^*}-1}(a_n),\delta/(10C))$ so that we have
$|F_{{m^*}-1}(z) - a_n|=R_n$.  By the choice of $\delta$ we see that
$h_{i_{m^*}}(z)$ is defined and since $|f_n(h_{i_{m^*}}(z))-a_n|=R_n$ we also
have $\delta/(10C) \leq |h_{i_{m^*}}(z)-k_n(a_n)|$ by the definition of $f_n$.
By the
choice of $C$ and property (III) we conclude $$\delta/(10C)
\leq |h_{i_{m^*}}(z)-k_n(a_n)| \leq \delta/10.$$
We note that (see Remark~\ref{gis1-1})
$f_n$ is one-to-one on $\tri(k_n(a_n),\delta)$.
By Lemma~\ref{domain} we have that $k_n$ is defined on
$\tri(a_n,2R_n)$ and Lemma~\ref{annlemma} shows that $k_n(A_n)$
contains the annulus $$B_n=Ann(k_n(a_n);r_n',R_n')$$ where
$R_n'=dist(k_n(C(a_n,R_n)),k_n(a_n))$ and
$r_n'=\max\{|k_n(z)-k_n(a_n)|:|z - a_n|=r_n\}$ satisfy
$R_n'/r_n' \geq
\frac{R_n}{9r_n}\to \infty.$  By construction
each $\delta/(10C) \leq R_n' \leq \delta/10$ (the lower bound follows
from the condition that $f_n(\tri(k_n(a_n),\delta/(10C))) \subset
\tri(a_n,R_n)$).  Each $B_n$ separates $A$ since $f_n(B_n) \subset
A_n$ and $f_n(A) \subset A$ (Lemma~\ref{forinv}).  Thus we see by
Lemma~\ref{Rntozero} that we have obtained the desired contradiction.
\end{proof}

\begin{proof}[Proof of Corollary~\ref{cor3}]

The corollary follows directly from Lemma~\ref{equivhyp} and
Theorem~\ref{better}.
\end{proof}


\begin{proof}[Proof of Corollary~\ref{main}]
Suppose $A$ contains two or more points.
Let $z \in A'$ with, say, $f(z)=z$ and let $w \in A
\setminus \{z\}$.  Since $f^n(w)$ are distinct points in $A$
(which tend to $z$) we see that
$A$ contains infinitely many points.  The result follows from
Corollary~\ref{cor3}.
\end{proof}

\begin{proof}[Proof of Corollary~\ref{cor}]
Since each $g_i$ for $i \in I$ is conformal, we
have that each $g_i'$ is never zero on $A$.   Since $A$ is compact
$\eta=\min_{i=1\dots N} \inf_{a \in A} |g_i'(a)|>0$ and so the result follows from Corollary~\ref{main}.
\end{proof}


\section{Examples}

\begin{example}
Let $G=\langle g_1, \dots, g_N \rangle$ where $g_i(z) = z^2 + c_i$ and
each $c_i$ is chosen as follows.  Fix $0<\epsilon<1/2$ and choose the
$c_i \in Ann(0;\epsilon^2,\epsilon - \epsilon^2) \subset
\tri(0,1/4)$ so that there exists a half plane $H$ with the origin on
its boundary which
contains $\cup_{i=1}^N \tri(c_i,\epsilon^2)$.  Choose a connected open
set $U$ which is contained $H$ such
that $U\supset \cup_{i=1}^N \tri(c_i,\epsilon^2)$ and $\overline{U}
\subset \tri(0,\epsilon)$.  Since each $g_i(\tri(0,\epsilon)) \subset U$
we see that each
$g_i(\overline{U}) \subset U.$ With this $U$ and $K=\cup_{i=1}^N
g_i(\overline{U})$ the IFS $G$ is conformal (even though the
generating maps are not globally conformal) and by Corollary~\ref{cor}
the resulting attractor set is uniformly perfect.
\end{example}

\begin{example}
Let $G=\langle g_i:i \in I \rangle$ where each $g_i(z) = z^2 + c_i$
and each $c_i$ is chosen as follows.  Fix $0<\epsilon<1/2$ and
$\rho>0$ such that $2\rho<\epsilon - 2 \epsilon^2$.  Choose $c_i \in
Ann(0;\rho+\epsilon^2,\epsilon - \epsilon^2 - \rho) \subset
\tri(0,1/4)$.  Hence for $U=\tri(0,\epsilon)$ and
$K=\overline{Ann(0;\rho,\epsilon -
\rho)}$ we have $g_i(U) = \tri(c_i,\epsilon^2) \subset K$.
Hence $A \subset K$ and we see
that for $\eta=2\rho$ and $a \in A$ each $|g_i'(a)|=2|a|\geq \eta$ for
$i \in I$.  By Corollary~\ref{cor3}, if $A$ has infinitely many points,
then $A$ is uniformly perfect.

If $A$ has three or more points, then $A$ has infinitely many points
thus $A$ is uniformly perfect.  This follows since $g_1(A) \subset
\tri(c_1,\epsilon^2)$ and so there exists a point $a \in A \cap
\tri(c_1,\epsilon^2)$ other than the attracting fixed point of $g_1$.
Hence the points of the form $g_1^n(a)$ are distinct since $g_1$ is
one-to-one on $\tri(c_1,\epsilon^2)$.  We see that $A$ has three or
more points if it is not the case that $G= \langle g_1, g_2 \rangle$
where, denoting the attracting fixed points of $g_i$ by $a_i$,
$a_2=\pm a_1$.  (If $G$ were such an IFS, then we would have
$A=\{a_1,a_2\}$.)  For example, one can check that if $c_1>0$
then $c_2 \neq c_1$ and $c_2 \neq c_1-1+\sqrt{1-4c_1}$ implies $A$ is
uniformly perfect.  Also, if three of the maps $g_i$ are distinct,
then $A$ is uniformly perfect.
\end{example}

\begin{example}
We construct a conformal IFS $G$ such that $A$ is not uniformly
perfect.  Choose sequences $a_n$ and $b_n$ of positive real numbers
such that $1>b_1>a_1>b_2>a_2>\dots$ where $a_n \downto 0, b_n \downto
0$ and $a_n/b_{n+1} \to \infty$.  Then the maps $g_n(z)=(b_n - a_n)z +
a_n$ are such that the conformal IFS $G=\langle g_n:n \geq 1 \rangle$
(where $U$ is an appropriately defined \nbhd of the interval $[0,1]$)
has an attractor set $A$ that is separated by each of the annuli
$Ann(0;b_{n+1},a_n)$.  Thus this attractor set is not uniformly
perfect.  We note that any attempt to apply Theorem~\ref{better}
breaks down since neither a choice of $\delta$ nor a choice of $C$ can
be found to satisfy its hypotheses.
\end{example}

\begin{example}\label{nonupfingen}
We construct an example of a finitely generated (polynomial) IFS which
has an attractor set which is not uniformly perfect.  Let
$f(z)=z^{23}$ and $g(z)=(z-1/2)^2 + 1/2$.  Let $I=[0,3/4],
I_{1/10}=\{z:|z-w| \le 1/10 \text{ for some } w \in I\}$ and let
$U=I_{1/10} \cup
\tri(1/2,0.37)$.  Since
$\overline{U} \subset \tri(0,0.9) \cap \tri(1/2,0.601)$
we see that $f(\overline{U}) \subset f(\tri(0,0.9)) \subset
\tri(0,1/10) \subset U$ and $g(\overline{U}) \subset
g(\tri(1/2,0.601))
\subset \tri(1/2,0.37)\subset U$.
Hence $G=\langle f,g \rangle$ is an IFS defined on $(U,K)$ where
$K=f(\overline{U}) \cup g(\overline{U})$.  Since $f(I) \cup g(I)
\subset I$ we see that $A \subset f(I) \cup g(I)$.  Let $V$ be the
open interval $((3/4)^{23},1/2)$ and note that $V \cap A =
\emptyset$.  We claim that $f^n(V) \cap A = \emptyset$ for all $n=1,
2, \dots$.

\begin{proof}
We first show that $f(V) \cap A = \emptyset$.  Suppose $a \in f(V)
\cap A$.  Since $A=f(A) \cup g(A)$ and
$f(V)\cap g(A) \subset f(I)\cap g(I) = \emptyset$ we conclude that $a
\in f(V) \cap f(A)$.  Hence there exists $a_1 \in A$ such that
$f(a_1)=a$.  Since $f$ is one-to-one on $I \supset A \cup V$ we must then
have $a_1 \in V$.  This contradiction shows that $f(V) \cap A =
\emptyset$.  One can then use this argument in an inductive manner to
conclude $f^n(V) \cap A =
\emptyset$ for all $n=1, 2, \dots$.

Since $f^n(V)=((3/4)^{23^{n+1}}, 2^{-23^n})$ we see that
$Ann(0;(3/4)^{23^{n+1}},2^{-23^n})$ separates $A$ and thus $A$ is
not uniformly perfect.  We note, however, that since
each map in $G$ is a polynomial and so has finite degree,
Lemma~\ref{Aperfect} implies $A$ is perfect since $A$ contains the
infinitely many points of the form $f^n(1/2)$.
\end{proof}

\begin{remark}
We note that in Example~\ref{nonupfingen} it can be shown that both
$f$ and $g$ are one-to-one on $A$, but not on any open set containing
$A$ and thus we may not apply Corollary~\ref{cor}.  To see that $g$ is
one-to-one on $A \subset f(I) \cup g(I)$ we note that $g(g(I)) \cap
g(f(I)) = \emptyset$ and both $g|_{g(I)}$ and $g|_{f(I)}$ are one-to-one.
\end{remark}
\end{example}


\begin{example}\label{finite}
We give an example that demonstrates that it is possible
for the attractor set of a (non conformal) IFS
which is generated by analytic maps to
have cardinality any finite number.

Fix a positive integer $n$.  Consider the complex polynomial
$f_0(z)= c+(z^n - c^n)^2$
which has a super-attracting fixed point at $c$, i.e., $f_0(c)=c$ and
$f_0'(c)=0$.  For
$\alpha=e^{2\pi i/n}$ we see that $f_0(\alpha^kc)=c$ for
$k=0,\dots,n-1$.  Hence $f_k(z)=\alpha^kf_0(z)$ has a super-attracting
fixed point at $\alpha^kc$ and maps each $\alpha^jc$ to $\alpha^kc$
for $j=0,\dots,n-1$.  It is elementary to verify that for small values
of $c$ and for a
suitable choice of $r>0$ we get that $f_k(\tri(0,r))\subset K$ for
$k=0,\dots,n-1$ where $K$ is a suitably chosen compact subset of ${\tri(0,r)}.$
Hence these functions generate an iterated function system on $(U,K)$
where $U=\tri(0,r)$.

One can easily see by~\eqref{invariance} that
$A=\{\alpha^kc:0 \leq k \leq n-1 \}$.
So for any $n$ there exists an attractor set $A$
such that $card(A)=n$.
\end{example}

\bibliographystyle{plain}


\end{document}